\documentclass{article}
\usepackage{amsmath,amssymb,amsfonts}
\usepackage{graphicx}
\usepackage{color}

\usepackage{pb-diagram}
\def\rank{\mathrm{rank}}

\def\inv{\mathrm{inv}}

\newtheorem{theorem}{Theorem}
\newtheorem{proposition}{Proposition}
\newtheorem{corollary}{Corollary}
\newtheorem{lemma}{Lemma}

\newtheorem{example}{Example}
\newtheorem{definition}{Definition}
\newenvironment{proof}[1][Proof]{\noindent\textit{#1.} }{\hfill$\Box$\medskip}

\title{Geometrization and Generalization of the Kowalevski top}
\author{Vladimir Dragovi\' c}

\date{}

\begin{document}

\maketitle

\medskip

\centerline{Mathematical Institute SANU}

\centerline{Kneza Mihaila 36, 11000 Belgrade, Serbia}

\smallskip

\centerline{Mathematical Physics Group, University of Lisbon}

\smallskip

\centerline{e-mail: {\tt vladad@mi.sanu.ac.rs}}

\

\begin{abstract}
\smallskip
A new view on the Kowalevski top and the Kowalevski integration
procedure is presented. For more than a century, the Kowalevski 1889
case, attracts full attention of a wide community as the highlight
of the classical theory of integrable systems. Despite hundreds of
papers on the subject, the Kowalevski integration is still
understood as a magic recipe, an unbelievable sequence of skilful
tricks, unexpected identities and smart changes of variables. The
novelty of our present approach is based on our four observations.
The first one is that the so-called fundamental Kowalevski equation
is an instance of a pencil equation of the theory of conics which
leads us to a new geometric interpretation of the Kowalevski
variables $w, x_1, x_2$ as the pencil parameter and the Darboux
coordinates, respectively. The second is observation of the key
algebraic property of the pencil equation which is followed by
introduction and study of a new class of {\bf discriminantly
separable polynomials}. All steps of the Kowalevski integration
procedure are now derived as easy and transparent logical
consequences of our theory of discriminantly separable polynomials.
The third observation connects the Kowalevski integration and the
pencil equation with the theory of multi-valued groups. The
Kowalevski change of variables is now recognized as an example of a
two-valued group operation and its action. The final observation is
surprising equivalence of the associativity of the two-valued group
operation and its action to $n=3$ case of the Great Poncelet Theorem
for pencils of conics.
\end{abstract}

\newpage

\tableofcontents

\newpage

\section{Introduction}\label{sec:intro}

\medskip
The goal of this paper is to give a new view on the Kowalevski top
and the Kowalevski integration procedure. For more than a century,
the Kowalevski 1889 case \cite{Kow}, attracts the full attention
of a wide community as the highlight of the classical theory of
integrable systems. Despite hundreds of papers on the subject, the
Kowalevski integration is still understood as a magic recipe, an
unbelievable sequence of skilful tricks, unexpected identities and
smart changes of variables (see for example \cite{Kow1},
\cite{Kot}, \cite{Gol}, \cite{Mlo}, \cite{App}, \cite{Del},
\cite{Jur}, \cite{Dub}, \cite{Aud}, \cite{BRST},\cite{VN},
\cite{DRV}, \cite{HM}, \cite{Kuz},\cite{Mar}, \cite{DG} and
references therein).

The novelty of this paper is based on our four observations. The
first one is that the so-called fundamental Kowalevski equation
(see \cite{Kow}, \cite{Kot}, \cite{Gol})
$$
Q(w, x_1, x_2)=0,
$$
is an instance of a pencil equation from the theory of conics.
This leads us to a new interpretation of the Kowalevski variables
$w, x_1, x_2$ as the pencil parameter and the Darboux coordinates
respectively. Origins and classical applications of the Darboux
coordinates can be found in Darboux's book \cite {Dar1} , while
some modern application can be found in \cite {Drag1}, \cite
{Drag2}.

The second is observation of the key algebraic property of the
pencil equation: {\it all three of its discriminants are expressed
as products of two polynomials in one variable each}:
$$
\aligned \mathcal D_w(Q)(x_1,x_2)&=f_1(x_1)f_2(x_2)\\
\mathcal D_{x_1}(Q)(w,x_2)&=f_3(w)f_2(x_2)\\
\mathcal D_{x_2}(Q)(w,x_1)&=f_1(x_1)f_3(w)
\endaligned
$$
This serves us as a motivation to introduce a new class of what we
call {\bf discriminantly separable polynomials}. We develop the
theory of such polynomials. All steps of the Kowalevski
integration now follow as easy and transparent logical
consequences of our theory of the discriminantly separable
polynomials.

The third observation connects the Kowalevski integration and the
pencil equation with the theory of multivalued groups. The theory of
multivalued groups started in the beginning of 1970's by Buchstaber
and Novikov (see \cite{BN}). It has been further developed by
Buchstaber and his collaborators in last forty years (see
\cite{Buc}, \cite{BV}, \cite{BR}). The Kowalevski change of
variables is now recognized as a case of two-valued group operation
$(\Gamma_2, \mathbb Z_2)$ and its action, where $\Gamma_2$ is an
elliptic curve and $\mathbb Z_2$ its subgroup.

Our final observation is  surprising equivalence of the
associativity condition for this two-valued group operation to a
case of the Great Poncelet Theorem for triangles. Well-known
mechanical interpretation of the Great Poncelet Theorem is connected
with integrable billiards, see for example \cite {DR}. The Great
Poncelet Theorem is the milestone of the theory of pencils of conics
and the whole classical projective geometry (see \cite{Pon}, and
also \cite{Ber}, \cite{DR}, \cite{DR2} and references therein), as
the Kowalevski top is the milestone of the classical integrable
systems. Now we manage to relate them closely. As a consequence, we
get a new connection between Great Poncelet Theorem and integrable
mechanical systems, this time from rigid- body dynamics.

The paper is organized as follows. The next Section \ref{sec:prel}
starts with a subsection devoted to the pencils of conics and the
Darboux coordinates. We derive the key property of the pencil
equation-discriminant separability. In the second subsection, we
formally introduce the class of discriminanlty separable
polynomials and systematically study this class.

In the Section \ref{sec:kow} we show how the Kowalevski case is
embedded into our more general framework. A new geometric
interpretation of the Kowalevski variables $(w, x_1, x_2)$ as the
pencil parameter and the Darboux coordinates is obtained.

In the Section \ref{sec:genintsys}  general  systems are defined,
related to the general equation of the pencil. The Kowalevski top
can be seen as a special subcase. The first integrals are studied.
Their properties are related to the properties of discriminantly
separable polynomials, obtained in Section \ref{sec:prel}. It was
done by use of what we call the Kotter trick (see \cite {Kot},
\cite{Gol}). The nature of this transformation is going to be
clarified in the last Section \ref{sec:dvg} through the theory of
multivalued groups. Then, we manage to generalize another Kotter's
transformation and this gives us a possibility to integrate the
general system defined at the beginning of this Section. We reduce
the problem to the functions $P_i, i=1,2,3$. The evolution of
those functions in terms of the theta-functions was obtained by
Kowalevski herself in \cite {Kow}. A modern account of the
theta-functions and their applications to nonlinear equations one
can find for example in \cite{Dub}.

The last Section \ref{sec:dvg} is devoted to two-valued groups and
their connection with the Kowalevski top and the Great Poncelet
Theorem. In order to make the text self-contained as much as
possible, we start the Section with brief introduction to the theory
of multivalued groups, following works of Buchstaber and his
co-workers. The main role is played by two-valued coset group
obtained from an elliptic curve $\Gamma_2$ and its subgroup $\mathbb
Z_2$. It appears that the Kowalevski change of variables has its
naturale expression through this two-valued group and its action.
These results complete the picture obtained before by Weil in
\cite{We} and Jurdjevic \cite{Jur}. Within this framework, we give
an explanation of the Kotter trick, as we promised in Section
\ref{sec:genintsys}. Finally, we show that the associativity
condition for the two-valued group $(\Gamma_2,\mathbb Z_2)$ is
equivalent to the famous Great Poncelet Theorem (\cite {Pon}) in its
basic  $n=3$ case.

\medskip

\section{Pencils of conics and discriminantly separable polynomials}\label{sec:prel}

\medskip

\subsection{Pencils of conics and the Darboux coordinates}
\medskip

Let us start with two conics $C_1$ and $C_2$ given by their
tangential equations:
\begin{equation}\label{eq:conics}
\aligned C_1:\,&
a_0w_1^2+a_2w_2^2+a_4w_3^2+2a_3w_2w_3+2a_5w_1w_3+2a_1w_1w_2=0;\\
 C_2:\,& w_2^2-4w_1w_3=0.
 \endaligned
 \end{equation}
We assume that conics $C_1$ and $C_2$ are in general position.
 Consider the pencil $C(s)$ of conics $C_1+sC_2$. The conics from the pencil share four common
 tangents. The
 coordinate equation of the conics of the pencil is:
 \begin{equation}\label{eq:defpolynom}
F(s,z_1,z_2,z_3):=\det M(s, z_1,z_2,z_3)=0,
\end{equation}
where $M$ is a bordered matrix of the form
\begin{equation}\label{eq:defmatrix}
M(s,z_1,z_2,z_3)=\left[\begin{array}{cccc} 0 & z_1 & z_2 & z_3\\
z_1 & a_0 & a_1 & a_5 - 2s\\
z_2 & a_1 & a_2 + s & a_3\\
z_3& a_5 - 2s & a_3 & a_4
\end{array}\right]
\end{equation}
\medskip
Then the point equation of the pencil of conics $C(s)$ is of the
form of the quadratic polynomial in $s$
\begin{equation}\label{eq:pencil}
F:=H+Ks+Ls^2=0
\end{equation}
where $H$, $K$ and $L$ are quadratic expressions in $(z_1, z_2,
z_3)$.

Following Darboux (see \cite{Dar1}), we introduce a new system of
coordinates in the plane. Given a plane with standard coordinates
$(z_1, z_2, z_3)$, we start from the given conic $C_2$. The conic
is given by the equation (\ref{eq:conics}) and it is rationally
parameterized by $(1, \ell, \ell^2)$. The tangent line to the
conic $C_2$ through the point with the parameter $\ell_0$ is given
by the equation
$$
t_{C_2}(\ell_0): z_1\ell_0^2-2z_2\ell_0+z_3=0.
$$
On the other hand, for a given point $P$ in the plane with
coordinates $P=(\hat z_1,\hat z_2,\hat z_3)$ there correspond two
solutions $x_1$ and $x_2$ of the equation quadratic in $\ell$:
\begin{equation}\label{eq:darboux}
\hat z_1\ell^2 - 2\hat z_2\ell + \hat z_3=0.
\end{equation}
Each solution corresponds to a tangent to the conic $C_2$ from the
point $P$. We will call the pair $(x_1, x_2)$ {\it the Darboux
coordinates} of the point $P$. One finds immediately converse
formulae
\begin{equation}\label{eq:viet}
\hat z_1=1,\quad \hat z_2= \frac{x_1+x_2}{2},\quad  \hat z_3 = x_1
x_2.
\end{equation}
\medskip
We change the variables in the polynomial $F$ from projective
coordinates $(z_1:z_2:z_3)$ to the Darboux coordinates according
to the formulae \ref{eq:viet}. In the new coordinates we get the
formulae:
\begin{equation}\label{eq:pencil1}
\aligned
H(x_1,x_2)=&(a_1^2-a_0a_2)x_1^2x_2^2+(a_0a_3-a_5a_1)x_1x_2(x_1+x_2)
\\&+(a_5^2-a_0a_4)(x_1^2+x_2^2)
+(2(a_5a_2-a_1a_3)+\frac{1}{2}(a_5^2-a_0a_4)x_1x_2\\&+(a_1a_4-a_3a_5))(x_1+x_2)+a_3^2-a_2a_4\\
K(x_1,x_2)=&-a_0x_1^2x_2^2+2a_1x_1x_2(x_1+x_2)-a_5(x_1^2+x_2^2)-4a_2x_1x_2\\
&+2a_3(x_1+x_2)-a_4\\
L(x_1,x_2)=&(x_1-x_2)^2.
\endaligned
\end{equation}
\medskip
We may notice for further references that
\begin{equation}\label{eq:difference}
(x_1-x_2)^2=4(z_1z_3-z_2^2).
\end{equation}
\medskip
Now, the polynomial
$$
F(s, x_1, x_2)=L(x_1,x_2)s^2 + K(x_1,x_2)s + H(x_1,x_2)
$$
is of the second degree in each of variables $s$, $x_1$ and $x_2$
and it is symmetric in $(x_1, x_2)$. It has one {\bf very
exceptional property}, as described in the next theorem.

For a polynomial $P(y_1,y_2,\dots, y_n)$ of variables
$(y_1,y_2,\dots, y_n)$ we will denote its discriminant with
respect to the variable $y_i$ by $\mathcal D_{y_i}(P)$ which is a
polynomial of the rest of the variables $(y_1,\dots,
y_{i-1},y_{i+1},\dots, y_n)$.
\medskip

\begin{theorem}\label{th:discrimsep}
\begin{itemize}
\item[(i)] There exists a polynomial $P=P(x)$ such that the
discriminant of the polynomial $F$  in $s$ as a polynomial in
variables $x_1$ and $x_2$ separates the variables:
\begin{equation}\label{eq:discsep1}
\mathcal D_{s}(F)(x_1,x_2)=P(x_1)P(x_2).
\end{equation}
\item[(ii)] There exists  a polynomial $J=J(s)$ such that the
discriminant of the polynomial $F$  in $x_2$ as a polynomial in
variables $x_1$ and $s$ separates the variables:
\begin{equation}\label{eq:discsep2}
\mathcal D_{x_2}(F)(s, x_1)=J(s)P(x_1).
\end{equation}
Due to the symmetry between $x_1$ and $x_2$ the last statement
remains valid after exchanging the places of $x_1$ and $x_2$.
\end{itemize}
\end{theorem}
\medskip
\begin{proof}
\begin{itemize}
\item[(i)] A general point belongs to two conics of a tangential
pencil. If a point belongs to only one conic, then it belongs to
one of the four common tangents of the pencil. At such a point,
this unique conic touches  one of the four common tangents. Thus,
the equation
\begin{equation}\label{eq:discsep3}
\mathcal D_s(F)(x_1,x_2)=0
\end{equation}
which represents the condition of annulation of the discriminant,
is the equation of the four common tangents. Thus, the equation
\ref{eq:discsep3} is equivalent to the system
$$
\aligned x_1&=c_1 \quad x_1=c_2 \quad x_1=c_3 \quad x_1=c_4\\
x_2&=c_1 \quad x_2=c_2 \quad x_2=c_3 \quad x_2=c_4
\endaligned
$$
where $c_i$ are parameters which correspond to the points of
contact of the four common tangents with the conic $C_2$. As  a
consequence, we get
$$
\mathcal D_s(F)(x_1,x_2)=P(x_1)P(x_2),
$$
where the polynomial $P$ is of the fourth degree and of the form
$$
P(x)=a(x-c_1)(x-c_2)(x-c_3)(x-c_4).
$$
This proves the first part of the theorem.

The second part of the Theorem follows from the following Lemma.
\medskip
\begin{lemma}\label{lemma:discsep}
Given  a polynomial $S=S(x,y,z)$ of the second degree in each of
its variables in the form:
$$
S(x,y,z)=A(y,z)x^2 + 2B(y,z)x + C(y,z).
$$
If there are polynomials $P_1$ and $P_2$ of the fourth degree such
that
\begin{equation}\label{eq:discsep4}
B(y,z)^2-A(y,z)C(y,z)=P_1(y)P_2(z),
\end{equation}
then there exists a polynomial $f$ such that
$$
\mathcal D_y S(x,z)=f(x)P_2(z), \quad D_z S(x,y)=f(x)P_1(y).
$$
\end{lemma}
\medskip
\begin{proof} To prove the Lemma, rewrite the equation
(\ref{eq:discsep4}) in the equivalent form
$$
(B+uA)^2-A(u^2A+2uB+C)=P_1(y)P_2(z).
$$
For a zero $y=y_0$ of the polynomial $P_1$, any zero of $S(u, y_0,
z)$ as a polynomial in $z$ is a double zero, according to the last
equation. Thus, $y_0$ is a zero of $\quad D_z S(x,y)$. Thus, the
polynomial $P_1$ is a factor of the polynomial $\quad D_z S(x,y)$.
Since the degree of the polynomial $P_1$ is four, then there
exists a polynomial $f$ in $x$ such that
$$
 D_z S(x,y)=f(x)P_1(y).
$$
The rest of the Lemma follows by double application of the same
arguments.
\end{proof}
\medskip
\item[(ii)] Now, the proof of the second part of the Theorem
\ref{th:discrimsep} follows by immediate application of the Lemma
\ref{lemma:discsep}
\end{itemize}
\end{proof}

\medskip
\begin{proposition}\label{prop:explicite}
\begin{itemize}
\item[(i)] The explicit formulae for the polynomials $P$ and $J$
are
\begin{equation}\label{eq:explicitplynom}
\aligned P(x)=& a_0x^4-4a_1x^3+(2a_5+4a_2)x^2 - 4 a_3x +a_4\\
J(s)=&-4s^3+4(a_5-a_2)s^2+(a_0a_4-a_5^2+4(a_5a_2-a_1a_3))s\\
&-a_3^2a_0+a_0a_4a_2+2a_1a_3a_5-a_4a_1^2-a_2a_5^2.
\endaligned
\end{equation}
\item[(ii)] If all the zeros of the polynomial $P$ are simple,
then the elliptic curves
$$
\aligned \Gamma_1: y^2&=P(x)\\
\Gamma_2: t^2&=J(s)
\endaligned
$$
are isomorphic and the later can be understood as Jacobian of the
former.
\end{itemize}
\end{proposition}
\medskip
\begin{proof} Instead of straightforward calculation, we are going
to consider a double-bordered determinant (see \cite{Dar1},
\cite{VD}, \cite{Hir}) obtained from the matrix $M$
(\ref{eq:defmatrix}):
\begin{equation}\label{eq:doublebord}
\hat M=\left|\begin{array}{ccccc} 0& 0 & z'_1 & z'_2 & z'_3\\
0& 0 & z_1 & z_2 & z_3\\
z'_1 & z_1 & a_0 & a_1 & a_5 - 2s\\
z'_2 & z_2 & a_1 & a_2 + s & a_3\\
z'_3 & z_3& a_5 - 2s & a_3 & a_4
\end{array}\right|
\end{equation}
We apply  the Jacobi identity and get
$$
\hat M_{11} \hat M_{22} - (\hat M_{12})^2= \hat M \hat M_{12,12}.
$$
Obviously, $\hat M_{12,12}$ is a polynomial  only in $s$ of the
third degree:
$$
\aligned \hat M_{12,12}=&-4s^3+4(a_5-a_2)s^2+((a_0a_4-a_5^2)
+4(a_5a_2-a_1a_3))s\\
&+a_0a_4a_2-a_3^2a_0+2a_1a_3a_5-a_4a_1^2-a_2a_5^2\\
=&J(s)
\endaligned
$$
Moreover, if we substitute
$$
\aligned z_1&=1,\quad  z_2= \frac{x_1+x_2}{2},\quad   z_3 = x_1
x_2\\
z'_1&=1,\quad  z'_2= \frac{x_1+x'_2}{2},\quad   z'_3 = x_1 x'_2
\endaligned
$$
we have
$$
\aligned \hat M&=P(x_1)\frac{(x_2-x'_2)^2}{4}\\
\hat M_{11}&=F(s,x_1,x_2)\\
 \hat M_{22}&=F(s,x_1,x'_2).
\endaligned
$$
If we denote
$$
F(s,x_1,x_2)=T(s,x_1)x_2^2+V(s,x_1)x_2 +W(s,x_1)
$$
then
$$
\hat M_{12}=Tx_2x'_2+V\frac{x_2+x'_2}{2}+W.
$$
From the last equations, after dividing by $(x_2-x'_2)^2$, we get
$$
V^2-4TW=J(s)P(x_1),
$$
and the proof of the first part of the Proposition is finished.

The second part follows by direct calculation of correspondence
between two elliptic curves, one of which is defined by a
polynomial of degree 3 and one by polynomial of degree 4.
\end{proof}
\medskip

\subsection{Discriminantly separable polynomials}
\medskip

We saw that a polynomial of three variables which defines a pencil
of conics has a very peculiar property: all three of its
discriminants are representable as products of two polynomials of
one variable each. These considerations motivate the following
definition.
\medskip
\begin{definition} For a polynomial $F(x_1,\dots,x_n)$ we say that it
is {\emph discriminantly separable} if there exist polynomials
$f_i(x_i)$ such that for every $i=1,\dots , n$
$$
\mathcal D_{x_i}F(x_1,\dots, \hat x_i, \dots, x_n)=\prod_{j\ne
i}f_j(x_j).
$$
It is {\emph symmetrically discriminantly separable} if
$$f_2=f_3=\dots = f_n,$$
while it is {\emph strongly discriminatly separable} if
$$
f_1=f_2=f_3=\dots = f_n.$$ It is {\emph weakly discriminantly
separable} if there exist polynomials $f^j_i(x_i)$ such that for
every $i=1,\dots , n$
$$
\mathcal D_{x_i}F(x_1,\dots, \hat x_i, \dots, x_n)=\prod_{j\ne
i}f^i_j(x_j).
$$
\end{definition}
\medskip
\begin{theorem} Given a polynomial $F(s,x_1,x_2)$ of the second
degree in each of the variables $s, x_1, x_2$ of the form
$$
F=s^2A(x_1,x_2)+2B(x_1,x_2)s + C(x_1,x_2).
$$
Denote by $T_{B^2-AC}$ a $5\times 5$ matrix such that
$$
(B^2-AC)(x_1,x_2)=\sum_{j=1}^5\sum_{i=1}^5T_{B^2-AC}^{ij}x_1^{i-1}x_2^{j-1}.
$$
Then, polynomial $F$ is discriminantly separable if and only if
$$
\rank \, T_{B^2-AC}=1.
$$
\end{theorem}
\medskip
\begin{proof} The proof follows from the Lemma \ref{lemma:discsep} and
the observation that a polynomial in two variables is equal to a
product of two polynomials in one variable if and only if its
matrix is equal to a tensor product of two vectors. The last
condition is equivalent to the condition on rank of the last
matrix to be equal to 1.
\end{proof}
\medskip
\begin{proposition}Given a polynomial $F(s,x_1,x_2)$ of the second
degree in each of the variables $s, x_1, x_2$. of the form
$$
F=s^2A(x_1)+2B(x_1,x_2)s + C(x_2),
$$
where $A$ depends only on $x_1$ and $C$ depends only on $x_2$.
Denote by $T_{B^2}$ a $5\times 5$ matrix such that
$$
(B^2)(x_1,x_2)=\sum_{i=1}^5T_{B^2}^{ij}x_1^{i-1}x_2^{j-1}.
$$
Then, polynomial $F$ is discriminantly separable if and only if
$$
\rank\, T_{B^2}=2.
$$
\end{proposition}
\medskip
\begin{proof} The proof follows from the observation of the proof
of the last theorem and a fact that a matrix of rank two is equal
to a sum of two matrices of rank one.
\end{proof}
\medskip

The last Proposition gives a method to construct nonsymmetric
discriminantly separable polynomials.
\medskip
\begin{lemma} Given an arbitrary quadratic polynomial
$$
F=s^2A+2Bs + C.
$$
Then, the square of its differential is equal to its discriminant
under the condition $F=0$:
$$
\left(\frac{dF}{ds}\right)^2=4(B^2-AC).
$$
\end{lemma}
\medskip
\begin{corollary} For an arbitrary discriminantly separable
polynomial $F(x_3,x_1,x_2)$ of the second degree in each of the
variables $x_3, x_1, x_2$, its differential is separable on the
surface $F(x_3,x_1,x_2)=0$:
$$
\frac{dF}{\sqrt{f_3(x_3)f_1(x_1)f_2(x_2)}}= \frac {dx_3}{\sqrt
{f_3(x_3)}}+\frac {dx_1}{\sqrt {f_1(x_1)}}+\frac {dx_2}{\sqrt
{f_2(x_2)}}.
$$
\end{corollary}
\medskip
The proof of the corollary is straightforward application of the
previous statements. This property of discriminantly separable
polynomials is fundamental in their role in the theory of
integrable systems. Observe that analogous statement is valid for
arbitrary discriminantly separable polynomials.

From the last Corollary, applied to a symmetric discriminatly
separable polynomial of the second degree, immediately follows a
variant of the Euler theorem.
\medskip
\begin{corollary} The condition $x_3= const $ defines a conic from the pencil
as an integral curve of the Euler equation:
$$
\frac {dx_1}{\sqrt {f_1(x_1)}}+\frac {dx_2}{\sqrt{ f_1(x_2)}}=0,
$$
where $f_1$ is general polynomial of degree $4$.
\end{corollary}
\medskip

\medskip
\begin{proposition}
All symmetric discriminantly separable polynomials $F(s, x_1,
x_2)$ of degree two in each variable with the leading coefficient
$$
L(x_1,x_2)=(x_1-x_2)^2
$$
are of the form
$$
F(s, x_1, x_2)=(x_1-x_2)^2s^2 + K(x_1,x_2)s + H(x_1,x_2)
$$
where $K$ and $H$ are done by the formulae (\ref{eq:pencil1}).
\end{proposition}
\medskip
The next Lemma gives a possibility to create new discriminantly
separable polynomials from a given one.
\medskip
\begin{lemma}\label{lemma:creat}
Given a discriminantly separable polynomial
$$
F(s, x_1, x_2):=A(x_1, x_2)s^2+2B(x_1,x_2)s+C(x_1,x_2)
$$
of the second degree in each variable.
\begin{itemize}
\item [(a)] Let $\alpha (x) $ be a linear transformation. Then
polynomial
$$
F_1(s, x_1, x_2):=F(s,\alpha (x_1) , x_2)
$$
is discriminantly separable. \item [(b)] The polynomial
$$
\hat F(s, x_1, x_2):=C(x_1, x_2)s^2+2B(x_1,x_2)s+A(x_1,x_2)
$$
is discriminantly separable.
\end{itemize}
\end{lemma}
\medskip
The transformation from $F$ to $\hat F$ described in the Lemma
\ref{lemma:creat} (b) maps a solution $s$ of the equation $F=0$ to
$1/s$. We will use the term {\it transposition} for such a
transformation from $F$ to $\hat F$. Thus, summarizing we get
\medskip
\begin{corollary}
Given a discriminantly separable polynomial
$$
F(s, x_1, x_2):=A(x_1, x_2)s^2+2B(x_1,x_2)s+C(x_1,x_2)
$$
of the second degree in each variable and three
fractionally-linear transformations $\alpha, \beta, \gamma $. Then
the polynomial
$$
F_1(s, x_1, x_2):=F(\gamma(s),\alpha (x_1) , \beta(x_2))
$$
is discriminantly separable.
\end{corollary}
\medskip
From the last Lemma we have a procedure to create non-symmetric
discriminantly separable polynomials from a given symmetric
discriminantly separable polynomial. The converse statement is
also true:
\medskip
\begin{proposition}
Given a discriminantly separable polynomial
$$
F(s, x_1, x_2):=A(x_1, x_2)s^2+2B(x_1,x_2)s+C(x_1,x_2)
$$
of the second degree in each variable. Suppose that a biquadratic
$F(s_0,x_1,x_2)$ is nondegenerate for some value $s=s_0$. Then
there exists a fractionally-linear transformations $\alpha$ such
that the polynomial
$$
F_1(s, x_1, x_2):=F(s,\alpha (x_1) , x_2)
$$
is symmetrically discriminantly separable.
\end{proposition}
\medskip
\begin{proof}
Let us fix an arbitrary value for $s$ such that $B(x_1,x_2)$ is a
nondegenerate biquadratic. Keeping $s$ fixed, we have a relation
$$
\frac {dx_1}{\sqrt {f_1(x_1)}}\pm\frac {dx_2}{\sqrt{ f_2(x_2)}}=0,
$$
where $f_1, f_2$ are two polynomials, each in one variable. For a
given $x_1$ there are two corresponding points $x_2$ and $\hat
x_2$. The last two are connected by the relation
$$
\frac {d\hat x_2}{\sqrt {f_2(\hat x_2)}}\pm\frac {dx_2}{\sqrt{
f_2(x_2)}}=0,
$$
where now denominator of both fraction is one and the same
polynomial, $f_2$. This means that there exists an elliptic
function $u$ of degree two and a shift $T$ on the elliptic curve
$y^2=f_2(x)$, such that $x_2$ and $\hat x_2$ are parameterized by
$$
x_2=u(z) \quad \hat x_2=u(z+T).
$$
From the relations
$$
B(x_1,x_2)=0 \quad B(x_1, \hat x_2)=0
$$
we see that both $y$ and $y^2$ are elliptic functions of degree at
most four which can be expressed through $x_2, \hat x_2$. Thus,
$y$ is an elliptic function of degree two. There is a
fractional-linear transformation which reduces $y$ to $u(z+T/2)$.
This concludes the proof of the Proposition.
\end{proof}
\medskip

\section{Geometric interpretation of the Kowalevski fundamental
equation}\label{sec:kow}
\medskip

The magic integration of the Kowalevski top is based on the
Kowalevski fundamental equation, see  \cite{Kot}, \cite{Gol}:
\begin{equation}\label{eq:kowrel1}
Q(w,x_1,x_2):=(x_1-x_2)^2w^2-2R(x_1,x_2)w-R_1(x_1,x_2)=0,
\end{equation}
where
\begin{equation}\label{eq:kowrel2}
\aligned R(x_1,x_2) =& -x_1^2x_2^2+6l_1x_1x_2+2lc(x_1+x_2)+c^2-k^2\\
R_1(x_1,x_2) =& -6l_1x_1^2x_2^2-(c^2-k^2)(x_1+x_2)^2-4clx_1x_2(x_1+x_2)\\
& + 6l_1(c^2-k^2)-4c^2l^2.
\endaligned
\end{equation}
If we replace in the equations (\ref{eq:pencil}) and
(\ref{eq:pencil1}) the following values for the coefficients:
\begin{equation}\label{eq:coeff}
\aligned a_0&=-2 \quad a_1=0 \quad a_5=0\\
a_2&=3l_1 \quad a_3=-2cl \quad a_4=2(c^2-k^2)
\endaligned
\end{equation}
and compare with (\ref{eq:kowrel1}) and (\ref{eq:kowrel2}), we get
the following
\medskip
\begin{theorem}\label{th:kowalevski}
The Kowalevski fundamental equation represents a point pencil of
conics given by their tangential equations
\begin{equation}\label{eq:conics1}
\aligned \hat C_1:\,&
-2w_1^2+3l_1w_2^2+2(c^2-k^2)w_3^2-4clw_2w_3=0;\\
 C_2:\,& w_2^2-4w_1w_3=0.
\endaligned
\end{equation}
The Kowalevski variables $w, x_1, x_2$ in this geometric settings
are the pencil parameter, and the Darboux coordinates with respect
to the conic $C_2$ respectively.
\end{theorem}
\medskip
The Kowalevski case corresponds to the general case under the
restrictions
$$
a_1=0 \quad a_5=0 \quad a_0=-2.
$$
The last of these three relations is just normalization condition,
provided  $a_0\ne 0$. The Kowalevski parameters $l_1, l, c$ are
calculated by the formulae
$$
l_1=\frac{a_2}{3} \quad
l=\pm\frac{1}{2}\sqrt{-a_4+\sqrt{a_4+4a_3^2}} \quad
c=\mp\frac{a_3}{\sqrt{-a_4+\sqrt{a_4+4a_3^2}}}
$$
provided  that $l$ and $c$ are requested to be real.
\medskip
Let us mention at the end of this Section, that in the original
paper \cite{Kow}, instead the relation (\ref{eq:kowrel1}),
Kowalevski used the equivalent one
$$
\hat
Q(s,x_1,x_2):=(x_1-x_2)^2(s-\frac{l_1}{2})^2-R(x_1,x_2)(s-\frac{l_1}{2})-\frac{R_1(x_1,x_2)}{4}=0.
$$
The equivalence is obtained by putting $w=2s-l_1$.
\medskip
\section{Generalized integrable system}\label{sec:genintsys}
\medskip
\subsection{Equations of motion and the first integrals}
\medskip
We are going to consider the following system of differential
equations on unknown functions $e_1, e_2, x_1, x_2, r, g$:
\begin{equation}\label{eq:diffsystem}
\aligned \frac {d e_1}{dt} &= -\alpha e_1\\
\frac {d e_2}{dt} &= \alpha e_2\\
\frac {d x_1}{dt} &= -\beta (rx_1+cg)\\
\frac {d x_2}{dt} &= \beta (rx_2+cg)\\
\frac {d r}{dt} &= -\beta(x_2- x_1)(x_1+x_2+a_1)-\frac{\alpha}{2r}(e_1-e_2)\\
\frac {d g}{dt} &= \frac{\beta}{2c}\left[(x_2-
x_1)(x_1x_2-a_5)+e_1x_2-e_2x_1\right]+\frac{(2r\beta-\alpha)}{2c^2g}\left(e_1x_2^2-e_2x_1^2\right)
\endaligned
\end{equation}
Here $\beta$ and $\alpha$ are given functions of $e_1, e_2, x_1,
x_2, r, g$. The choice of their form defines different systems.
The Kowalevski top is equivalent to the above system for
$$
a_1=0 \quad a_5=0,
$$
with the choice
\begin{equation}\label{eq:kowab}
\alpha = i r \quad \beta =\frac{i}{2}.
\end{equation}
We will assume in what follows that $a_1$ and $a_5$ are general.
Beside the last choice for $\alpha$ and $\beta$, there  are many
others choices which also provide polynomial  vector fields, such
as (A) $\alpha=kr^2 \quad \beta=\frac{k}{2}r$, (B) $\alpha=krg
\quad \beta=k_1g$, (C) $\alpha=kr^2g \quad\beta=k_1g$. Interesting
cases  satisfy the system (\ref{eq:measure}) from Proposition
(\ref{prop:measure}).
\medskip
\begin{proposition}\label{prop:integrals}
The system (\ref{eq:diffsystem}) has the following first integrals
\begin{equation}\label{eq:integrals}
\aligned k^2&=e_1 \cdot e_2\\
a_0a_2&=e_1 + e_2 - (x_1+x_2)^2-2a_1(x_1+x_2) - r^2\\
-\frac{a_0a_3}{2}&=-x_2e_1-x_1e_2+x_1x_2(x_1+x_2)+\frac{a_5}{2}(x_1+x_2)
+
a_1x_1x_2-rg\\
\frac{a_0a_4}{4}&=x_2^2e_1+x_1^2e_2-x_1^2x_2^2-a_5x_1x_2-g^2
\endaligned
\end{equation}
\end{proposition}
\medskip
One can rewrite the last relations in the following form
\begin{equation}\label{eq:integral2}
\aligned k^2&=e_1 \cdot e_2\\
r^2&=e_1 + e_2  + \hat E(x_1,x_2)\\
rg&=-x_2e_1-x_1e_2 + \hat F(x_1, x_2)\\
g^2&=x_2^2e_1+x_1^2e_2 + \hat G(x_1,x_2),
\endaligned
\end{equation}
where
\begin{equation}\label{eq:EFG}
\aligned
\hat E(x_1,x_2)&=-a_0a_2- K (x_1+x_2)^2-2a_1(x_1+x_2)\\
\hat F(x_1,x_2)&=\frac{a_0a_3}{2}+ K
x_1x_2(x_1+x_2)+\frac{a_5}{2}(x_1+x_2)
+a_1x_1x_2   \\
\hat G(x_1,x_2)&=-\frac{a_0a_4}{4}-K x_1^2x_2^2-a_5x_1x_2,
\endaligned
\end{equation}
with
$$
K=1.
$$
\medskip
\begin{lemma}\label{lemma:EFG}
 If the polynomials $\hat E, \hat F, \hat G$ are defined by the
equation (\ref{eq:EFG}) then the polynomial
$$
P(x_1):=\hat E(x_1,x_2)x_1^2+2\hat F(x_1,x_2)x_1+\hat G(x_1,x_2)
$$
depends only on $x_1$.
\end{lemma}
\medskip
\begin{proposition} Given three polynomials $\hat E(x_1,x_2),
\hat F(x_1,x_2), \hat G(x_1, x_2)$ of the second degree in each
variable such that
\begin{itemize}
\item[(1)] Polynomials $P, Q$ defined by
\begin{equation}\aligned
P(x_1):=&\hat E(x_1,x_2)x_1^2+2\hat F(x_1,x_2)x_1+\hat G(x_1,x_2)\\
Q(x_2):=&\hat E(x_1,x_2)x_2^2+2\hat F(x_1,x_2)x_2+\hat G(x_1,x_2)
\endaligned
\end{equation}
depend only on one variable each. \item[(2)] Polynomials
$R(x_1,x_2)$ and $R_1(x_1,x_2)$ defined by
\begin{equation}\aligned
R(x_1,x_2):=&\hat E(x_1,x_2)x_1x_2+\hat F(x_1,x_2)(x_1+x_2)+\hat G(x_1,x_2)\\
R_1(x_1,x_2):=&\hat E(x_1,x_2)\hat G(x_1,x_2) - \hat F^2(x_1,x_2)
\endaligned
\end{equation}
are of the second degree in each variables. \end{itemize}
 Then:
 \begin{itemize}
 \item[(a)] The polynomials $\hat E(x_1,x_2),
\hat F(x_1,x_2), \hat G(x_1, x_2)$  are symmetric in $x_1, x_2$.
\item[(b)] The polynomial
$$
F(s, x_1,x_2)=(x_1-x_2)^2s^2-2R(x_1,x_2)s - R_1(x_1,x_2)
$$
is discriminantly separable. \item[(c)] The most general form of
the polynomials $\hat E,\hat F,\hat G$ is given in the equation
(\ref{eq:EFG}), with $K$ arbitrary.
\item[(d)] For $K=1$  the
polynomial $P$ is the one given in the Proposition
\ref{prop:explicite}.
\end{itemize}
\end{proposition}
\medskip
\begin{proof}The proof follows by straightforward calculation with
application of the Lemma \ref{lemma:discsep}.
\end{proof}
\medskip

If the coefficient $K$ is nonzero we may normalize it to be equal
to one. Under this assumption, the equations (\ref{eq:EFG}) with
$K=1$ are general. The case $K=0$ is going to be analyzed
separately in one of the following sections.

 From the equations (\ref{eq:integral2}) we get the following
\medskip
\begin{corollary}
The relation is satisfied
\begin{equation}\label{eq:identity1}
e_2P(x_1)+e_1P(x_2)-H(x_1,x_2)+k^2(x_1-x_2)^2=0.
\end{equation}
where $P$ is the polynomial defined in the Lemma \ref{lemma:EFG}.
\end{corollary}
\medskip
\begin{corollary} The differentials of $x_1$ and $ x_2$ may be
written in the form
\begin{equation}\label{eq:dx1dx2}
\aligned \frac{dx_1}{dt}&=-\beta\sqrt{P(x_1)+e_1(x_1-x_2)^2}\\
\frac{dx_2}{dt}&=\beta\sqrt{P(x_2)+e_2(x_1-x_2)^2}.
\endaligned
\end{equation}
\end{corollary}
\medskip
The proof follows from the equations (\ref{eq:integral2}) and
Lemma \ref{lemma:EFG}.
\medskip
Now, we apply what we are going to call {\it the Kotter trick}:
\medskip
\begin{equation}\label{eq:kotter1}
\left[\sqrt{e_1}\frac{\sqrt{P(x_2)}}{x_1-x_2}\pm
\sqrt{e_2}\frac{\sqrt{P(x_1)}}{x_1-x_2}\right]^2=(w_1\pm k)(w_2
\mp k),
\end{equation}
where $w_1, w_2$ are solutions of the quadratic equation
\begin{equation}\label{eq:kotter2}
F(s, x_1,x_2)=(x_1-x_2)^2s^2-2R(x_1,x_2)s - R_1(x_1,x_2).
\end{equation}
\medskip
The Kotter trick appeared in \cite{Kot} quite mysteriously.
Further explanation done by Golubev sixty years later seems to be
even trickier, see \cite{Gol} and much less clear. In the last
section of this paper, see Proposition \ref{prop:commdiagram}, we
provide a new interpretation of this transformation as a commuting
diagram of morphisms of double-valued group. Should we hope that
our explanation is more transparent then previous ones, since new
sixty years passed in meantime?

From the last relations, following Kotter, one gets
$$
\aligned
\left(\frac{dx_1}{\sqrt{P(x_1)}dt}\right)^2&=\beta^2\left(1+\frac{(x_1-x_2)^4e_1P(x_2)}{P(x_1)P(x_2)(x_1-x_2)^2}\right)\\
&=\beta^2\left(1+\frac{(\sqrt{(w_1-k)(w_2+k)}+\sqrt{(w_1+k)(w_2-k)})^2}{(w_1-w_2)^2}\right)\\
\left(\frac{dx_2}{\sqrt{P(x_2)}dt}\right)^2&=\beta^2\left(1+\frac{(x_1-x_2)^4e_2P(x_1)}{P(x_1)P(x_2)(x_1-x_2)^2}\right)\\
&=\beta^2\left(1+\frac{(\sqrt{(w_1-k)(w_2+k)}-\sqrt{(w_1+k)(w_2-k)})^2}{(w_1-w_2)^2}\right).
\endaligned
$$
Next, we get
\begin{equation}\label{eq:dxdt}
\aligned
\frac{dx_1}{\sqrt{P(x_1)}dt}&=-\beta\left(\frac{\sqrt{(w_1-k)(w_1+k)}+\sqrt{(w_2+k)(w_2-k)}}{(w_1-w_2)}\right)\\
\frac{dx_2}{\sqrt{P(x_2)}dt}&=-\beta\left(\frac{\sqrt{(w_1-k)(w_1+k)}-\sqrt{(w_2+k)(w_2-k)}}{(w_1-w_2)}\right)
\endaligned
\end{equation}
\medskip
Now we apply the discriminant separability property of the
polynomial $F$:
\begin{equation}\label{eq:kowchange}
\aligned
\frac{dx_1}{\sqrt{P(x_1)}}+\frac{dx_2}{\sqrt{P(x_2)}}&=\frac{dw_1}{\sqrt{J(w_1)}}\\
\frac{dx_1}{\sqrt{P(x_1)}}-\frac{dx_2}{\sqrt{P(x_2)}}&=\frac{dw_2}{\sqrt{J(w_2)}}
\endaligned
\end{equation}
\medskip
We will refer to the last relations as {\it the Kowalevski change
of variables}. The nature of these relations has been studied by
Jurdjevic (see \cite{Jur}) following Weil (\cite {We}). We are
going to develop further these efforts in the Section
\ref{sec:dvg} where we are going to show that the Kowalevski
change of variables is the infinitesimal version of a double
valued group operation and its action.
\medskip
From the relations \ref{eq:kowchange} and \ref{eq:dxdt} we finally
get:
\begin{equation}\label{eq:kowvareq}
\aligned
\frac{dw_1}{\sqrt{\Phi(w_1)}}+\frac{dw_2}{\sqrt{\Phi(w_2)}}&=0\\
\frac{w_1\,dw_1}{\sqrt{\Phi(w_1)}}+\frac{w_2\,dw_2}{\sqrt{\Phi(w_2)}}&=2\beta\,dt,
\endaligned
\end{equation}
where
$$
\Phi(w)=J(w)(w-k)(w+k),
$$
is the polynomial of  fifth degree. Thus, the equations
(\ref{eq:kowvareq}) represent the Abel-Jacobi map of the genus $2$
curve
$$
y^2=\Phi(w).
$$
\medskip
\subsection{Generalized Kotter transformation}
\medskip

In order to integrate the dynamics on the Jacobian of the
hyper-elliptic curve $y^2=\Phi(w)$ we are going to generalize
classical Kotter transformation. In this section we will assume
the normalization condition
$$
a_0=-2.
$$
\begin{proposition}
For the polynomial $F(s, x_1, x_2)$ there exist polynomials
$A_0(s)$, $f(s)$, $A(s, x_1,x_2)$, $B(s, x_1, x_2)$ such that  the
following identity
\begin{equation}
F(s,x_1,x_2)\cdot A_0(s)=A^2(s,x_1,x_2)+f(s)\cdot B(s,x_1,x_2),
\end{equation}
is satisfied. The polynomials are defined by the formulae:
$$
\aligned A(s,x_1,x_2)&=A_0(s)(x_1x_2-s)+B_0(s)(x_1+x_2)+M_0(s)\\
A_0(s)&=a_1^2-a_0a_2-sa_0\\
B_0(s)&=\frac{1}{2}(a_0a_3-a_5a_1+2sa_1)\\
M_0(s)&=a_5a_2-a_1a_3+s(a_1^2+a_5)\\
B(s,x_1,x_2)&=(x_1+x_2)^2+2a_1(x_1+x_2)-2s-2a_2\\
f(s)&=2s^3+2(a_2-a_5)s^2+\left(2(a_1a_3-a_5a_2)+a_4+\frac{a_5^2}{2}\right)s
+f_0\\
f_0&=a_4a_2-a_3^2-a_1a_3a_5+\frac{a_4a_1^2+a_2a_5^2}{2}.
\endaligned
$$
\end{proposition}
\medskip
For $a_5=a_1=0$ the previous identity has been obtained in
\cite{Kot}. Following Kotter's idea, consider the identity
$$
F(s)=F(u)+(s-u)F'(u)+(s-u)^2.
$$
From the last two identities we get a quadratic equation in $s-u$
$$
(s-u)^2(x_1-x_2)^2-2(s-u)(R(x_1,x_2)-u(x_1-x_2)) +f(u)B
+(x_1-x_2)^2A^2.
$$
\medskip
\begin{corollary}
\begin{itemize}
\item[(a)] The solutions of the last equation satisfy the identity
in $u$:
$$
(s_1-u)(s_2-u)=\frac{A^2}{(x_1-x_2)^2}+f(u)\frac{B}{(x_1-x_2)^2}.
$$
\item[(b)] Denote $m_1, m_2, m_3$ the zeros of the polynomial $f$,
and
$$P_{i}=\sqrt{(s_1-m_i)(s_2-m_i)}, \quad i=1, 2, 3.$$
Then
\begin{equation}\label{eq:Pi1}
P_i=\frac{1}{x_1-x_2}\left(\sqrt{A_0(m_i)}x_1x_2+\frac{B_0(m_i)}{\sqrt{A_0(m_i)}}
+m_i(m_i-a_5-2a_2)-2a_5-a_1a_3\right),\quad i=1,2,3.
\end{equation}
\end{itemize}
\end{corollary}
\medskip
Now we introduce more convenient notation
$$
\aligned n_i&=m_i+a_1^2+2a_2, \quad i=1, 2, 3;\\
X&=\frac{x_1x_2+(2a_1^2+a_5+2a_2)+\frac{a_1}{2}(x_1-x_2)}{x_1-x_2},\\
Y&=\frac{1}{x_1-x_2},\\
Z&=\frac{(a_1^3+2a_2a_1+2a_5a_1+2a_3)(x_1+x_2)-2(a_1^2+2a_2)(a_1^2+a_5)}{x_1-x_2}.
\endaligned
$$
\medskip
\begin{lemma}
The quantities $X,Y,Z$ satisfy the system of linear equations
\begin{equation}\label{eq:linsys}
\aligned
X-n_1Y+\frac{1}{2n_1}Z&=\frac{P_1}{\sqrt{n_1}}\\
X-n_2Y+\frac{1}{2n_2}Z&=\frac{P_2}{\sqrt{n_2}}\\
X-n_3Y+\frac{1}{2n_3}Z&=\frac{P_3}{\sqrt{n_3}}.
\endaligned
\end{equation}
\end{lemma}
\medskip
Denote
$$
\hat f(x)=f(x-a_1^2-2a_2).
$$
One can easily solve the previous linear system and get
\begin{lemma}
The solutions of the system (\ref{eq:linsys}) are
$$
\aligned Y&=-\left(\frac{P_1\sqrt{n_1}}{\hat
f'(n_1)}+\frac{P_2\sqrt{n_2}}{\hat
f'(n_2)}+\frac{P_3\sqrt{n_3}}{\hat f'(n_3)}\right)\\
Z&=2n_1n_2n_3\left(\frac{P_1}{\sqrt{n_1}\hat
f'(n_1)}+\frac{P_2}{\sqrt{n_2}\hat
f'(n_2)}+\frac{P_3}{\sqrt{n_3}\hat f'(n_3)}\right)
\endaligned
$$
\end{lemma}
\medskip
The expression in terms of theta functions for
$P_i=\sqrt{(s_1-m_i)(s_2-m_i)}$ for $i=1, 2, 3$ can be obtained
from \cite{Kow} paragraph 7.
\medskip

\subsection{Interpretation of the equations of motion}
\medskip

\centerline{\bf Rigid-body coordinates}
\medskip

We are going to present briefly the interpretation of the equations
of motion (\ref{eq:diffsystem}) in the standard rigid-body
coordinates $p, q, r, \gamma, \gamma', \gamma''$, where:
$$
\aligned
 e_1&=x_1^2+c(\gamma +i\gamma')\\
 e_2&=x_2^2+c(\gamma -i\gamma')\\
 p&=\frac{x_1+x_2}{2}\\
q&=\frac{x_1-x_2}{2i}.
\endaligned
$$
From the last four equations of the system (\ref{eq:diffsystem})
we get
\medskip
\begin{equation}\label{eq:diffsystemclass1}
\aligned \dot p=& -i \beta  rq\\
\dot q=&i \beta rp \\
\dot r=&2\beta iq(2p+a_1)-\frac{i\alpha}{r}(2pq+c\gamma')\\
\dot \gamma''=& -\frac{\beta}{c}(qia_5+2i c\gamma
q-2ic\gamma' p)\\
&+\frac{2r\beta-\alpha}{c^2\gamma''}(ic\gamma'(p^2-q^2)-2icpq\gamma)
\endaligned
\end{equation}
while the equations for $\dot \gamma, \dot \gamma'$ can easily be
obtained from the first two equations of the system
(\ref{eq:diffsystem}):
\medskip
$$
\aligned \dot \gamma=&\frac{\alpha}{2c}(x_2^2-x_1^2)-i\alpha
\gamma' +\frac{-x_1\dot x_1-x_2\dot x_2}{c}\\
\dot \gamma'=&\frac{\alpha}{2c}(-x_2^2-x_1^2)-i\alpha \gamma
+\frac{-x_1\dot x_1+x_2\dot x_2}{c}.
\endaligned
$$
\medskip
Finally, we get
\begin{equation}\label{eq:diffsystemclass2}
\aligned \dot \gamma =&\frac {2i(2\beta r - \alpha)}{c}pq -i\alpha
\gamma' +2i\beta \gamma'' q\\
\dot \gamma' =&-\frac {2i(2\beta r - \alpha)}{c}(p^2-q^2) +i\alpha
\gamma -2i\beta \gamma'' q
\endaligned
\end{equation}
\medskip
\begin{proposition}\label{prop:measure}
The system (\ref{eq:diffsystemclass1}, \ref{eq:diffsystemclass2})
preserves the standard measure if and only if
\begin{equation}\label{eq:measure}
\aligned &A_0\alpha + A_1\alpha_p+A_2\alpha_q + A_3\alpha_r +
A_4\alpha_{\gamma} + A_5\alpha_{\gamma'} +A_6\alpha_{\gamma''} +\\
&B_0\beta + B_1\beta_p + B_2\beta_q+ B_3\beta_r+
B_4\beta_{\gamma}+ B_5\beta_{\gamma'} +B_6\beta_{\gamma''}=0,
\endaligned
\end{equation}
where
$$
\aligned
A_0&=r^2\gamma'p^2+c^2\gamma''^2\gamma'-2r^2pq\gamma+2c\gamma''^2pq-r^2\gamma'q^2\\
A_1&=0\\
A_2&=0\\
A_3&=-2c\gamma''^2rpq-c^2\gamma''^2r\gamma'\\
A_4&=-2pqr^2\gamma''^2-\gamma'r^2c\gamma''^2\\
A_5&=-2r^2\gamma''^2q^2+gr^2c\gamma''^2+2r^2\gamma''^2p^2\\
A_6&=-r^2\gamma''\gamma'p^2+2r^2\gamma''pq\gamma+r^2\gamma''\gamma'q^2\\
B_0&=-2r^3\gamma'p^2+2r^3\gamma'q^2+4r^3pq\gamma\\
B_1&=-cr^3q\gamma''^2\\
B_2&=cr^3p\gamma''^2\\
B_3&=4qr^2c\gamma''^2p+2qr^2c\gamma''^2a_1 \\
B_4&=2\gamma''^3qr^2c+4pqr^3\gamma''^2 \\
B_5&= -4r^3\gamma''^2p^2-2\gamma''^3qr^2c+4r^3\gamma''^2q^2\\
B_6&=
-r^2\gamma''^2qa_5-2r^3\gamma''\gamma'q^2-2r^2\gamma''^2c\gamma
q+2r^3\gamma''\gamma'p^2+2r^2\gamma''^2c\gamma'p-4r^3\gamma''pq\gamma
\endaligned
$$
\end{proposition}
\medskip
\begin{example} From the Kowalevski case, there is a pair $ \alpha = ir, \beta = i/2$ which satisfies the system (\ref{eq:diffsystemclass2})
written above. We give two more pairs:
$$\alpha_1= 2r(p^2+q^2) \quad \beta_1 = p^2+q^2,$$
and
$$\alpha_2=r\gamma'' \quad \beta_2=0.$$
 Moreover, any linear combination of the pairs
$(\alpha,\beta)$, $(\alpha_1,\beta_1)$ and $(\alpha_2,\beta_2)$ also
gives a solution of the system (\ref{eq:diffsystemclass2}) and
provides a system with invariant standard measure.
\end{example}

\medskip

\centerline{\bf Elastic deformations}

\medskip

Jurdjevic considered a deformation of the Kowalevski case associated
to a Kirchhoff elastic problem, see \cite{Jur}. The systems are
defined by the Hamiltonians
$$
H=M_1^2+M_2^2 + 2M_3^2 +\gamma_1
$$
where deformed Poisson structures $\{\cdot, \cdot\}_{\tau}$ are
defined by
$$
\{M_i, M_j\}_{\tau}=\epsilon_{ijk}M_k,\quad \{M_i,
\gamma_j\}_{\tau}=\epsilon_{ijk}\gamma_k,\quad \{\gamma_i,
\gamma_j\}_{\tau}=\tau \epsilon_{ijk}M_k,
$$
where the deformation parameter takes values $\tau = 0, 1, -1$. The
classical Kowalevski case corresponds to the case $\tau=0$.

Denote
$$
\aligned e_1&=x_1^2-(\gamma_1 + i\gamma_2) +\tau\\
e_2&=x_2^2-(\gamma_1 - i\gamma_2) +\tau,
\endaligned
$$
where
$$
x_{1,2}=\frac{M_1\pm i M_2}{2}.
$$
The integrals of motion
$$
\aligned I_1&=e_1 e_2\\
I_2&=H\\
I_3&=\gamma_1 M_1 + \gamma_2 M_2 + \gamma_3 M_3\\
I_4&=\gamma_1^2 + \gamma_2^2+\gamma_3^2 + \tau (M_1^2+M_2^2+M_3^2)
\endaligned
$$
may be rewritten in the form (\ref{eq:integral2})
$$
\aligned k^2&=I_1=e_1 \cdot e_2\\
M_3^2&=e_1 + e_2  + \hat E(x_1,x_2)\\
M_3\gamma_3&=-x_2e_1-x_1e_2 + \hat F(x_1, x_2)\\
\gamma_3^2&=x_2^2e_1+x_1^2e_2 + \hat G(x_1,x_2),
\endaligned
$$
where
$$
\aligned \hat G(x_1,x_2)&=-x_1^2x_2^2-2\tau x_1x_2 -2\tau
(I_1-\tau)+\tau^2-I_2\\
\hat F(x_1,x_2)&=(x_1x_2+\tau)(x_1+x_2)+I_3\\
\hat E(x_1,x_2)&=-(x_1+x_2)^2+2(I_1-\tau).
\endaligned
$$
\begin{proposition}
Corresponding pencil of conics is determined by equations
$$
a_1=0,\, a_5=2\tau,\,
a_2=\frac{2(\tau-I_1)}{a_0},\,a_3=2\frac{I_3}{a_0},\,a_4=\frac{8\tau(I_1-\tau)+4(I_2-\tau^2)}{a_0}
$$
where $a_0$ is arbitrary.
\end{proposition}
\medskip
\section{Two-valued groups, Kowalevski
equation and Poncelet Porism}\label{sec:dvg}
\medskip
\subsection{Multivalued groups: defining notions}
\medskip

The structure  of multivalued groups was introduced by Buchstaber
and Novikov in 1971 (see \cite{BN}) in their study of
characteristic classes of vector bundles, and it has been studied
by Buchstaber and his collaborators since then (see \cite{Buc} and
references therein).

Following \cite{Buc}, we give the definition of an n-valued group
on $X$ as a map:
$$
\aligned
 &m:\, X\times X \rightarrow (X)^n\\
 &m(x,y)=x*y=[z_1,\dots, z_n],
 \endaligned
 $$
where $(X)^n$ denotes the symmetric $n$-th power of $X$ and $z_i$
coordinates therein.

{\it Associativity} is the condition of equality of two $n^2$-sets
$$
\aligned &[x*(y*z)_1,\dots, x*(y*z)_n]\\
&[(x*y)_1*z,\dots, (x*y)_n*z]
\endaligned
$$
for all triplets $(x,y,z)\in X^3$.

An element $e\in X$ is {\it a unit} if
$$
e*x=x*e=[x,\dots,x],
$$
for all $x\in X$.

A map $\inv: X\rightarrow X$ is {\it an inverse} if it satisfies
$$
e\in \inv(x)*x, \quad e\in x*\inv(x),
$$
for all $x\in X$.

Following Buchstaber, we say that $m$ defines {\it an $n$-valued
group structure} $(X, m, e, \inv)$ if it is associative, with a
unit and an inverse.

An $n$-valued group $X$ acts on the set $Y$ if there is a mapping
$$
\aligned &\phi:\, X\times Y \rightarrow (Y)^n\\
&\phi (x,y)=x\circ y,
\endaligned
$$
such that the two $n^2$-multisubsets of $Y$
$$
x_1\circ (x_2\circ y) \quad (x_1*x_2)\circ y
$$
are equal for all $x_1, x_2\in X, y\in Y$. It is additionally
required that
$$
e\circ y=[y,\dots, y]
$$
for all $y\in Y$.
\medskip
\begin{example}[A two-valued group structure on $\mathbb {Z}_+$, \cite{BV}]
Let us consider the set of nonnegative integers $\mathbb Z_+$ and
define a mapping
$$
\aligned &m:\, \mathbb Z_+ \times \mathbb Z_+ \rightarrow (\mathbb
Z_+)^2,\\
&m(x,y)=[x+y,|x-y|].
\endaligned
$$
This mapping provides a structure of a two-valued group on
$\mathbb Z_+$ with the unit $e=0$ and the inverse equal to the
identity $\inv(x)=x$.

In \cite{BV} sequence of two-valued mappings associated with the
Poncelet porism was identified as the algebraic representation of
this 2-valued group. Moreover, the algebraic action of this group
on $\mathbb {CP}^1$ was studied and it was shown that in the
irreducible case all such actions are generated by Euler-Chasles
correspondences.
\end{example}

In the sequel, we are going to show that there is another 2-valued
group and its action on $\mathbb {CP}^1$ which is even more
closely related to the Euler-Chasles correspondence and to the
Great Poncelet Theorem, and which is at the same time intimately
related to the Kowalevski fundamental equation and to the
Kowalevski change of variables.

However, we will start our approach with a simple example.
\medskip

\subsection{The simplest case: 2-valued group $p_2$}
\medskip
Among the basic examples of multivalued groups, there are
$n$-valued additive group structures on $\mathbb C$. For $n=2$,
this is a two-valued group $p_2$ defined by the relation
\begin{equation}\label{eq:p2}
\aligned &m_2:\, \mathbb C \times \mathbb C \rightarrow (\mathbb
C)^2\\
&\, x *_2 y =[(\sqrt{x}+\sqrt{y})^2, (\sqrt{x}-\sqrt{y})^2]
\endaligned
\end{equation}
\medskip
The product $x *_2 y$ corresponds to the roots in $z$ of the
polynomial equation
$$
p_2(z, x, y)=0,
$$
where
$$
p_2(z, x, y)= (x+y+z)^2-4(xy+yz+zx).
$$
Our starting point in this section is the following
\medskip
\begin{lemma}\label{lemma:p2}
The polynomial $p_2(z, x, y)$ is discriminantly separable. The
discriminants satisfy relations
$$
\mathcal D_z(p_2)(x,y)=P(x)P(y) \quad \mathcal
D_x(p_2)(y,z)=P(y)P(z) \quad D_y(p_2)(x,z)=P(x)P(z),
$$
where
$$
P(x)=2x.
$$
\end{lemma}
\medskip
The polynomial $p_2$ as  discriminantly separable, generates a
case of generalized Kowalevski system of differential equations,
but this time  with $K=0$. The system is defined by
\begin{equation}\label{eq:sysp2}
\hat E=0 \quad \hat F=1 \quad \hat G=0,
\end{equation}
\medskip
and the equations of motion have the form
\begin{equation}\label{eq:sysp3}
\aligned \frac {d e_1}{dt} &= -\alpha e_1\\
\frac {d e_2}{dt} &= \alpha e_2\\
\frac {d x_1}{dt} &= -\beta (rx_1+cg)\\
\frac {d x_2}{dt} &= \beta (rx_2+cg)\\
\frac {d r}{dt} &= -\frac{\alpha}{2r}(e_1-e_2)\\
\frac {d g}{dt} &= 2\beta
c+\frac{(2r\beta-\alpha)}{2c^2g}\left(e_1x_2^2-e_2x_1^2\right)
\endaligned
\end{equation}
\medskip
In the standard rigid-body coordinates with $\alpha = i r$,
$\beta=i/2$ the last two equations become
$$
\dot r = 2pq + c \gamma' \quad \dot \gamma '' = ic.
$$
\medskip

\medskip
\begin{lemma}\label{lemma:p22}
The integrals of the system defined by the equations
(\ref{eq:sysp2}) are
$$
\aligned k^2&=e_1e_2\\
r^2&= e_1+e_2\\
crg&= 1-x_1e_2-x_2e_1\\
c^2g^2&=x_2^2e_1+x_1^2e_2
\endaligned
$$
\end{lemma}
\medskip
From the last Lemma \ref{lemma:p22} we get the relation
$$
2e_1x_2+2e_2x_1-1+k^2(x_1-x_2)^2=0.
$$
Now, together with the first integral relation from the Lemma
\ref{lemma:p22}, similar as in the Kowalevski case, we get
\begin{equation}\label{eq:sysp22}
\left[\sqrt{e_1}\frac{\sqrt{2x_2}}{x_1-x_2}\pm
\sqrt{e_2}\frac{\sqrt{2x_1}}{x_1-x_2}\right]^2=(w_1\pm k)(w_2 \mp
k),
\end{equation}
where $w_1, w_2$ are solutions of the quadratic equation
\begin{equation}\label{eq:sysp23}
F_2(w, x_1, x_2):=(x_1-x_2)^2w^2-2(x_1+x_2)w+1=0.
\end{equation}
The polynomial $F_2$ is obtained by transposition from the
polynomial $p_2$ and, thus, it is discriminantly separable:
$$
\mathcal D_x(F_2)(y,z)=P(y)\varphi (z),
$$
where
$$
\varphi (z)=z^3.
$$
Following lines of integration, we finally come to
\medskip
\begin{proposition}
The system of differential equations defined by \ref{eq:sysp2} is
integrated  to through the solutions of the system
\begin{equation}\label{eq:sysp22}
\aligned \frac {ds_1}{s_1\sqrt{\Phi_1(s_1)}} +\frac
{ds_2}{s_2\sqrt{\Phi_1(s_2)}}&=0\\
\frac {ds_1}{\sqrt{\Phi_1(s_1)}} +\frac
{ds_2}{\sqrt{\Phi_1(s_2)}}&=\frac{i}{2}dt,
\endaligned
\end{equation}
where
$$\Phi(s)=s(s-e_4)(s-e_5)$$
is the polynomial of degree 3.
\end{proposition}
\medskip
Similar systems appeared in a  slightly different context in the
works of Appel'rot, Mlodzeevskii, Delone in their study of
degenerations of the Kowalevski top (see \cite {App},
\cite{Mlo},\cite {Del}). In particular, we may construct {\it
Delone-type} solutions of the last system:
$$
s_1=0, \quad s_2=\wp\left(\frac{i}{4}(t-t_0)\right).
$$
\medskip

We can also consider integrable perturbation of the previous
integrable system, defined by:
\begin{equation}\label{eq:persys}
\aligned \hat E&=k_1-2a_1(x_1+x_2)\\
\hat F&=k_2+\frac{a_5}{2}(x_1+x_2)+a_1x_1x_2\\
\hat G&=k_3-a_5x_1x_2.
\endaligned
\end{equation}
\medskip
The equations of motion have the form
\begin{equation}\label{eq:sysp4}
\aligned \frac {d e_1}{dt} &= -\alpha e_1\\
\frac {d e_2}{dt} &= \alpha e_2\\
\frac {d x_1}{dt} &= -\beta (rx_1+cg)\\
\frac {d x_2}{dt} &= \beta (rx_2+cg)\\
\frac {d r}{dt} &= -\frac{\alpha}{2r}(e_1-e_2)-\frac{a_1}{2}\beta (x_2-x_1)\\
\frac {d g}{dt} &= 2\beta
c+\frac{(2r\beta-\alpha)}{2c^2g}\left(e_1x_2^2-e_2x_1^2\right)+\frac{a_5}{2}c\beta
(x_2-x_1)
\endaligned
\end{equation}
\medskip
In the standard rigid-body coordinates with $\alpha = i r$,
$\beta=i/2$ the last two equations become
$$
\aligned \dot r =& 2pq + c \gamma'+ \frac{a_1}{2} q \\
\dot \gamma '' = &ic(1+ i\frac{ a_5}{2} q).
\endaligned
$$
\medskip

Corresponding polynomial
$$
F(s, x_1, x_2)= (x_1-x_2)^2s^2 - 2R(x_1,x_2)s - R_1(x_1,x_2)
$$
where
$$
R(x_1,x_2)=\hat Ex_1x_2 + \hat F(x_1+x_2)+\hat G,\quad
R_1(x_1,x_2)=\hat E \hat F-\hat G^2,
$$
is discriminantly separable and
$$
\mathcal D_{x_1}(s, x_2)=\varphi(s)P(x_2),
$$
where
$$
\aligned \varphi(s)&=(2s-a_5)(2a_1+a_5s-2s^2)\\
P(x)&=2x(2a_1x^2-a_5x-2).
\endaligned
$$
\medskip

\subsection{2-valued group structure on $\mathbb {CP}^1$,
the Kowalevski fundamental equation and Poncelet porism}

\medskip

Now we pass to the general case. We are going to show that the
general pencil equation represents an action of a two valued group
structure. Recognition of this structure enables us to give to 'the
mysterious Kowalevski change of variables' a final algebro-geometric
expression and explanation, developing further the ideas of Weil and
Jurdjevic (see \cite{We}, \cite{Jur}). Amazingly, the associativity
condition for this action from geometric point of view is nothing
else than the Great Poncelet Theorem for a triangle.

As we have already mentioned, the general pencil equation
$$
F(s,x_1, x_2)=0
$$
is connected with two isomorphic elliptic curves
$$
\aligned \Gamma_1: y^2&=P(x)\\
\Gamma_2: t^2&=J(s)
\endaligned
$$
where the polynomials $P, J$ of degree four and three respectively
are defined by the equations (\ref{eq:explicitplynom}). Suppose
that the cubic one $\Gamma_2$ is rewritten in the canonical form
$$
\Gamma_2: t^2=J'(s)=4s^3-g_2s-g_3.
$$
Moreover, denote by $\psi:\, \Gamma_2\rightarrow \Gamma_1$ a
birational morphism between the curves induced by a
fractional-linear transformation $\hat \psi$ which maps three
zeros of $J'$ and $\infty$ to the four zeros of the polynomial
$P$.

The curve $\Gamma_2$ as a cubic curve has the group structure.
Together with its subgroup $\mathbb {Z}_2$ it defines the standard
two-valued group structure of coset type on $\mathbb {CP}^1$ (see
\cite {BR}, \cite{Buc}):
\begin{equation}\label{eq:G2Z2}
s_1 *_c s_2 =
\left[-s_1-s_2+\left(\frac{t_1-t_2}{2(s_1-s_2)}\right)^2,-s_1-s_2+\left(\frac{t_1+t_2}{2(s_1-s_2)}\right)^2\right],
\end{equation}
where $t_i=J'(s_i), i=1,2.$
\medskip
\begin{theorem}\label{th:G2Z21}
The general pencil equation after fractional-linear
transformations
$$
F(s, \hat \psi^{-1}(x_1), \hat \psi^{-1}(x_2))=0
$$
defines the two valued coset group structure $(\Gamma_2, \mathbb
Z_2)$ defined by the relation (\ref{eq:G2Z2}).
\end{theorem}
\medskip
\begin{proof}
After the fractional-linear transformations, the pencil equation
obtains the form
$$
F_1(s,x, y)=T(s,x)y^2+V(s,x)y + W(s,x),
$$
where
$$
\aligned T(s,x)&=-4s^2+4sx-s^2\\
V(s,x)&=4sx^2+2s^2x-2xg_2-g_2s-4g_3\\
W(s,x)&=-s^2x^2-g_2xs-4xg_3-2g_3s-\frac{g_2^2}{4}.
\endaligned
$$
We apply now a linear change of variables $\gamma$ on $s$:
$$
m=\gamma (s):= \frac{s}{2}
$$
and get
$$
F_2(m,x,y)=F_1(2m, x,y).
$$
Denote by $P=(m,n)$ and $M=(x,u)$  two arbitrary points on the
curve $\Gamma_2$, which means
$$
\aligned n^2&=4m^3-g_2m-g_3\\
u^2&=4x^3-g_2x-g_3.
\endaligned
$$
We want to find points $N_1=(y_1,v_1)$ and $N_2=(y_2,v_2)$ on
$\Gamma_2$ which correspond by $F_2$ to $P$ and $M$. These points
are
$$
\aligned y_1&=\frac{-V(s,x)+4nu}{2T(s,x)} \quad v_1=-\frac{2xT(s,
y_1)+V(s,y_1)}{4n}\\
y_2&=\frac{-V(s,x)-4nu}{2T(s,x)} \quad v_2=-\frac{2xT(s,
y_2)+V(s,y_2)}{4n}
\endaligned
$$
By trivial algebraic transformations
$$
\aligned y_1&=\frac{-4mx^2-4xm^2+xg_2+mg_2+2g_3+2nu}{-4(x-m)^2}\\
&=\frac{-4mx(x+m)+x^3 +m^3 - x^3 + xg_2+ g_3 -m^3 + mg_2+g_3+2nu}{-4(x-m)^2}\\
&=-x-m+\left(\frac{u-n}{2(x-m)}\right)^2
\endaligned
$$
we get the first part of the operation of the two-valued group
$(\Gamma_2, \mathbb Z_2)$ defined by the relation (\ref{eq:G2Z2}).
Applying similar transformations to $y_2$ we get the second part of
the relation (\ref{eq:G2Z2}) as well. This ends the proof of the
Theorem.
\end{proof}
\medskip
The Kowalevski change of variables (see equations
(\ref{eq:kowchange})) is infinitesimal of the correspondence which
maps a pair of points $(M_1,M_2)$ from the curve $\Gamma_1$ to a
pair of points $(S_1, S_2)$ of the curve $\Gamma_2$. One view to
this correspondence has been given in \cite {Jur} following Weil
\cite{We}. In our approach, there is a geometric view to this
mapping as the correspondence which maps {\it two tangents to the
conic $C$ to the pair of conics from the pencil which contain the
intersection point of the two lines}.

If we apply fractional-linear transformations to transform the curve
$\Gamma_1$ into the curve $\Gamma_2$, then the above correspondence
is nothing else then the two-valued group operation $*_c$ on
$(\Gamma_2, \mathbb Z_2)$.

\medskip
\begin{theorem}\label{th:G2Z22}
The Kowalevski change of variables is equivalent to infinitesimal of
the action of the two valued coset group $(\Gamma_2,\mathbb Z_2)$ on
$\Gamma_1$. Up to the fractional-linear transformation, it is
equivalent to the operation of the two valued group
$(\Gamma_2,\mathbb Z_2)$.
\end{theorem}
\medskip
Now, the Kotter trick from the Section \ref{sec:genintsys} (see
the equations (\ref{eq:kotter1}, \ref{eq:kotter2}) can be
presented as a commutative diagram.
\medskip
\begin{proposition}\label{prop:commdiagram}
The Kotter transformation defined by the equations
(\ref{eq:kotter1}, \ref{eq:kotter2}) makes the following diagram
commutative:
$$
\begin{diagram}
\node{\mathbb C^4}
 \arrow{e,t}{i_{\Gamma_1} \times i_{\Gamma_1}\times m}
 \arrow{s,l}{i_{\Gamma_1} \times i_{\Gamma_1} \times id \times id}
 \arrow{se,t}{i_a\times i_a \times m}
\node{\Gamma_1 \times \Gamma_1\times \mathbb C}
 \arrow{s,r}{p_1 \times p_1 \times id}
 \arrow{e,t}{\psi^{-1}\times  \psi^{-1}\times id}
\node{\Gamma_2 \times \Gamma_2\times \mathbb C}
 \arrow{ssw,b}{p_1 \times p_1 \times id}
  \\
\node{\Gamma_1 \times \Gamma_1\times \mathbb C \times \mathbb C}
 \arrow{s,l}{\varphi_1\times\varphi_2}
\node{\mathbb {CP}^1 \times \mathbb {CP}^1\times \mathbb C}%
 \arrow{s,l}{\hat \psi^{-1}\times \hat \psi^{-1}\times id}
  \\
\node{\mathbb C \times \mathbb C}
 \arrow{s,l}{m_2}
\node{\mathbb {CP}^1 \times \mathbb {CP}^1\times \mathbb C}
 \arrow{s,l}{m_c \times \tau_c}
  \\
\node{\mathbb {CP}^2}%
\node{\mathbb {CP}^2 \times \mathbb C/\sim}
 \arrow{w,t}{f}
\end{diagram}
$$

The mappings are defined as follows
$$
\aligned i_{\Gamma_1}&:\,x\mapsto (x,\sqrt{P(x)})\\
m&:\,(x,y) \mapsto  x\cdot y\\
i_a&:\, x \mapsto (x,1)\\
p_1&:\, (x, y)\mapsto x \\
m_c&:\, (x, y)\mapsto x*_c y\\
\tau_c&:\, x\mapsto (\sqrt{x},-\sqrt{x})\\
 \varphi_1&:\, (x_1,x_2, e_1, e_2)\mapsto
\sqrt{e_1}\frac{\sqrt{P(x_2)}}{x_1-x_2}\\
\varphi_2&:\, (x_1,x_2, e_1, e_2)\mapsto
\sqrt{e_2}\frac{\sqrt{P(x_1)}}{x_1-x_2}\\
f&:\, ((s_1,s_2,1),(k,-k))\mapsto
[(\gamma^{-1}(s_1)+k)(\gamma^{-1}(s_2)-k),(\gamma^{-1}(s_2)+k)(\gamma^{-1}(s_1)-k)]
\endaligned
$$
\end{proposition}
\medskip
From the Proposition \ref{prop:commdiagram} we see that the
two-valued group plays an important role in the Kowalevski system
and its generalizations.

 Putting together the
geometric meaning of the pencil equation and algebraic structure of
the two valued group we come to the connection with the Great
Poncelet Theorem (\cite {Pon}, see also \cite{Ber}, \cite {DR} and
\cite {DR2}). For the reader's sake we are going to formulate the
Great Poncelet Theorem for triangles in the form we are going to use
below.
\medskip
\begin{theorem}[Great Poncelet Theorem for triangles \cite {Pon}]
Given four conics $C_1, C_2, C_3, C$ from a pencil and three lines
$a_1, a_2, a_3$, tangents to the conic $C$ such that $a_1, a_2$
intersect on $C_1$, $a_2, a_3$ intersect on $C_2$ and $a_2, a_3$
intersect on $C_3$. Moreover, we suppose that the tangents to the
conics $C_1, C_2, C_3$ at the intersection points  are not
concurrent. Given $b_1, b_2$ tangents to the conic $C$ which
intersect at $C_1$. Then there exists $b_3$, tangent to the conic
$C$ such that the triplet $(b_1,b_2,b_3)$ satisfies all conditions
as $(a_1,a_2,a_3)$.
\end{theorem}
\medskip
Now, we are going back to the associativity condition for the
action of the double-valued group $(\Gamma_2,\mathbb Z_2)$.
\medskip
\begin{theorem}\label{th:G2Z23}
Associativity conditions for the group structure of the two-valued
coset group $(\Gamma_2,\mathbb Z_2)$ and for its action on
$\Gamma_1$ are equivalent to the great Poncelet theorem for a
triangle.
\end{theorem}
\medskip
\begin{proof}
Denote by $P$ and $Q$ two arbitrary elements of the two-valued group
$(\Gamma_2,\mathbb Z_2)$ and $M$ an arbitrary point on the curve
$\Gamma_1$. Let
$$
Q*P=[P_1, P_2]
$$
and
$$
P\circ M=[N_1, N_2].
$$
Associativity means the equality of the two quadruples:
$$
[Q\circ N_1, Q\circ N_2]=[P_1\circ M, P_2\circ M].
$$
Let us consider previous situation from geometric point of view.
Recall the geometric meaning  of the equation of a pencil of
conics
$$
F(s, x_1, x_2)=0.
$$
Variables $x_1$ and $x_2$ denote the Darboux coordinates of two
tangents to the conic $C_2$ which intersect at the conic $C_s$
with the pencil parameter equal to $s$.

Denote by $C_P$ and $C_Q$ the conics from the pencil which
correspond to the elements $P, Q$, and by $l_M, l_{N_1}, l_{N_2}$
the tangents to the conic $C_2$ which correspond to the points $M,
N_1, N_2$ of the curve $\Gamma_1$. Then, $l_{N_1}$ and $l_{N_2}$
are the two lines tangent to $C_2$ which intersect $l_M$ at the
conic $C_P$.

Moreover, if we denote
$$
Q\circ N_1=[N_3, N_4], \quad Q\circ N_2=[N_5, N_6]
$$
then corresponding lines $l_{N_3}, l_{N_4}, l_{N_5}, l_{N_6}$,
tangent to the conic $C_2$ satisfy the conditions: the pairs of
lines $(l_{N_1}, l_{N_3})$, $(l_{N_1}, l_{N_4})$, $(l_{N_2},
l_{N_5})$, $(l_{N_2}, l_{N_6})$ all intersect at the conic $C_Q$.

Now, associativity of the action is equivalent to the existence of
{\it a pair of conics $(C_{P_1}, C_{P_2})$} such that $(l_{M},
l_{N_3})$ and $(l_{M}, l_{N_6})$ intersect at the conic $C_{P_1}$,
while  $(l_{M}, l_{N_5})$ and $(l_{M}, l_{N_4})$ intersect at the
conic $C_{P_2}$, see the Fig. \ref{fig:assoc}.

\begin{figure}[h]\label{fig:assoc}
\centering
\includegraphics[width=8.8cm,height=6.1cm]{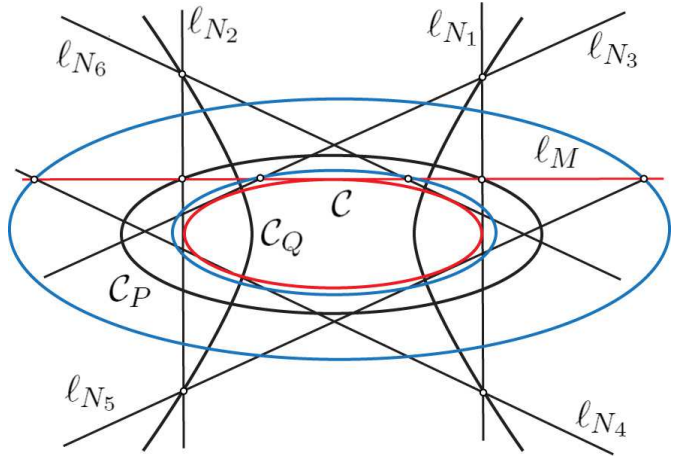}
\caption{Associativity condition and Poncelet theorem}
\end{figure}

Consider the intersection of the lines $(l_{M}, l_{N_3})$. Choose
the conic from the pencil which contains the intersection point,
such that the tangent to this conic at the intersection point is
not concurrent with the tangents to the conics $C_P$ and $C_Q$ at
the intersection points $(l_{M}, l_{N_1})$ and $(l_{N_1},
l_{N_3})$ respectively. Denote the conic $C_{P_1}$. Then by
applying Great Poncelet Theorem for triangle (see the Theorem
above, \cite {Pon},see also \cite {Ber}, \cite {DR}, \cite{DR2}),
one of the lines $l_{N_5}$ and $l_{N_6}$, say the last one,
intersects $L_{M}$  at the conic $C_{P_1}$. The tangent to this
conic at the intersection point is not concurrent with the
tangents to the conics $C_P$ and $C_Q$ at the intersection points
$(l_{M}, l_{N_2})$ and $(l_{N_2}, l_{N_6})$ respectively.

In the same way, by considering intersection of the lines $(l_{M},
l_{N_4})$ we come to the conic $(C_{P_2})$ from the pencil, which,
by Great Poncelet Theorem contains intersections of $(l_{M},
l_{N_4})$ and $(l_{M}, l_{N_5})$.

Since the result of the operation in the double-valued group
between elements $P, Q$  doesn't depend on the choice of the point
$M$ to which the action is applied, the conics $C_{P_2}$ and
$C_{P_1}$ in the previous construction should not depend of the
choice of the line $l_M$. This independence is equivalent to {\bf
the poristic} nature of the Poncelet Theorem. This demonstrates
the equivalence between the associativity condition and the Great
Poncelet Theorem for a triangle.
\end{proof}
\medskip
\medskip
From the last two theorems we get finally
\medskip

{\bf Conclusion} {\it Geometric settings for the Kowalevski change
of variables is the Great Poncelet Theorem for a triangle.}
\medskip

\subsection*{Acknowledgements}

The author is grateful to Borislav Gaji\' c and Katarina Kuki\' c
for helpful remarks. The research was partially supported by the
Serbian Ministry of Science and Technology, Project {\it Geometry
and Topology of Manifolds and Integrable Dynamical Systems}.  A part
of the paper has been written during a visit to the IHES. The author
uses the opportunity to thank the IHES for hospitality and
outstanding working conditions.

\newpage\thispagestyle{empty}
\vspace*{20mm}

\end{document}